\documentclass[11pt]{article}
\usepackage{amsmath,amsthm,amsfonts,latexsym,amssymb,enumerate,color}
\usepackage{graphicx}

\def\CC{\mathbb C}
\def\DD{\mathbb D}

\def\G{\mathcal G}

\def\TT{\mathbb T}

\def\ker{\mathop{\rm ker}\nolimits}

\newtheorem{thm}{Theorem}[section]
\newtheorem{prop}[thm]{Proposition}
\newtheorem{defn}[thm]{Definition}
\newtheorem{lem}[thm]{Lemma}
\newtheorem{cor}[thm]{Corollary}
\newtheorem{rem}[thm]{Remark}

\numberwithin{equation}{section}

\def\beginpf{\begin{proof}}
\def\endpf{\end{proof}}
\def\beq{\begin{equation}}
\def\eeq{\end{equation}}

\def\ol{\overline}

\def\sab{S_{a,b}}
\def\sigab{\Sigma_{a,b}}
\def\sigba{\Sigma_{b,a}}
\def\k{{\rm ker}}
\def\kerab{\k_{a,b}}
\def\kerba{\k_{b,a}}
\def\H2b{\overline{H^2_0}}

\def\ta{\tilde a}
\def\tb{\tilde b}
\def\dq{\"} 
\def\Hol{\mathop{\rm Hol}\nolimits}

\begin{document}

\title{Paired operators and paired kernels}

\author{M.~Cristina C\^amara,\thanks{
Center for Mathematical Analysis, Geometry and Dynamical Systems,
Instituto Superior T\'ecnico, Universidade de Lisboa, 
Av. Rovisco Pais, 1049-001 Lisboa, Portugal.
 \tt ccamara@math.ist.utl.pt} 
 \and Andr\'e Guimar\~{a}es,\thanks{Instituto Superior T\'ecnico, Universidade de Lisboa, 
 Av. Rovisco Pais, 1049-001 Lisboa, Portugal.
 \tt andre.guimaraes@tecnico.ulisboa.pt}
 \and  Jonathan R.~Partington\thanks{School of Mathematics, University of Leeds, Leeds LS2~9JT, U.K. {\tt j.r.partington@leeds.ac.uk}}\  
}

\date{}

\maketitle

\begin{abstract}
This paper is concerned with paired operators in the context of the Lebesgue
Hilbert space on the unit circle and its subspace, the Hardy space.
By considering when such operators commute,
generalizations of the Brown--Halmos results for Toeplitz operators are derived.
Further, the kernels of such operators are described, giving results on invariant and
nearly-invariant subspaces, together with a generalization of Coburn's theorem on
Toeplitz kernels.
\end{abstract}



\section{Introduction}

Let $X$ be a Hilbert space and  $P,Q$
complementary orthogonal projections on $X$. For $A,B$ bounded linear operators on $X$  we define the paired operators
 $S_{A,B}=AP+BQ$ and  $\Sigma_{A,B}=PA+QB$ (see \cite{MP}). These
first appeared in the  context of singular integral equations, and have been the object of
renewed interest because of the recently-established
connection with problems in diffraction theory \cite{speck}, dual truncated Toepliz operators
\cite{CKP}, and the characterization of block Toeplitz operators whose kernels present a scalar-type structure \cite[Cor.~3.4]{CP20}.

Here we work in the particular context of $L^2(\TT)$, where $\TT$ is the unit circle, and
$P$ is $P^+$, the orthogonal projection onto the Hardy space $H^2$ of the unit disc and $Q$ is 
$P^-$, the complementary projection onto $(H^2)^\perp=\H2b$.
The operators $A$ and $B$ are multiplication operators $A=M_a$ and $B=M_b$
with $a,b \in L^\infty(\TT)$.
This is the most common setting where these operators are studied, and also that with the most applications. They have been studied usually under an extra condition
that either $a$ or $b$ is invertible in $L^\infty(\TT)$ or,
at most, has zeros of ``integer order'' at a finite
set of points of $\TT$ (see \cite[Sec.~VI.8]{MP} for a more precise definition).

Paired operators are also closely connected with Toeplitz operators: in fact, they are dilations of
Toeplitz operators.
Moreover, it can be shown that if $b^{\pm 1} \in L^\infty(\TT)$, then the paired operator
$aP^+ + bP^-$ on $L^2(\TT)$ is equivalent after extension (as defined in \cite{BTsk}) to
the Toeplitz operator $T_{a/b}$ on $H^2$ (see \cite{C}), and an analogous
result can be obtained for dual Toeplitz operators  if $a^{\pm 1} \in L^\infty(\TT)$.
Thus it is natural to ask whether various well known properties of Toeplitz operators have analogues for paired operators.

In this paper we explore new connections of paired operators
with Toeplitz and Hankel operators. We study their
algebraic properties, such as when two such operators commute,
providing a generalization of the Brown--Halmos results for Toeplitz operators.
Then we consider the kernels of paired operators: this theme
leads off in two directions, one of which concerns invariant and nearly-invariant subspaces, 
while in the other direction we arrive at a generalization of Coburn's theorem 
on Toeplitz kernels.

\section{Paired operators: basic properties}

From now on we write $L^p$ to denote $L^p(\TT)$ for $p=2$ and $p=\infty$.

Let $a,b \in L^\infty$. The {\em paired operator} $\sab: L^2 \to L^2$ is
defined by 
\beq
\sab f = aP^+f + b P^- f, \qquad (f \in L^2).
\eeq
We call $\sigab: L^2 \to L^2$, defined by
\beq
\sigab f = P^+af + P^- bf, \qquad (f \in L^2),
\eeq
a {\em transposed paired operator}. We say that $\sab$ and $\sigab$ are {\em dual} to each other,
while $\sab$ and $\sigba$ are said to be {\em transposed} to each other (see \cite{MP}).

 \begin{defn}\label{def:ndg}
 We  say that a pair $\{a,b\}$ is {\em nondegenerate\/}
if $a, b \in L^\infty$ are nonzero a.e.\ on $\TT$, and   $a-b$ is too.
\end{defn}
Nondegeneracy will be a common hypothesis in our work, although sometimes
the weaker condition that $a -b \ne 0$ will be sufficient.

Note that we identify $a$ with $M_a$, the operator of multiplication by $a$ whenever the context
is clear. In particular, we write $P^+a+P^-b$ instead of
$P^+aI+P^-bI$.

Using the orthogonal decomposition $L^2=H^2 \oplus \overline{H^2_0}$  we may write down
the matrix representation of $S_{a,b}$. For convenience, write $P^\pm a=a_\pm$ and similarly for $b$.
Now since for $h_+ \in H^2$ we have $S_{a,b}h_+=ah_+$, the matrix form
of $\sab$ is a matrix of Hankel and Toeplitz operators:
\[
S_{a,b}=\begin{pmatrix}
a_+ + P^+a_- & P^+ b_ + \\
P^- a_- & P^-b_++b_- 
\end{pmatrix}: H^2 \oplus \overline{H^2_0} \to H^2 \oplus \overline{H^2_0}.
\]
Note that initially this matrix is only densely defined, but since
it corresponds to a bounded operator it can be extended
to the whole of $L^2$.\\

These various operators are related as follows.

\begin{prop}
(i) $\sab^*=\Sigma_{\bar a,\bar b}$, and\\
(ii) $\overline{\sab f}=zS_{\bar b,\bar a} \bar z \bar f
= z \Sigma_{b,a}^* \bar z \bar f$
for all $f \in L^2$.
\end{prop}

\beginpf 
(i) is clear; (ii) follows from the relations $\overline {P^\pm f }=zP^\mp \bar z \bar f$.
\endpf

\begin{prop}
\label{prop:2.2}
Let $M= \max ( \|a\|_\infty, \|b\|_\infty)$. Then
\beq\label{eq:2.2}
M \le \|\sab\| \le \min(\sqrt2 M, \|a\|_\infty + \|b\|_\infty).
\eeq
\end{prop}

\beginpf
Restricting to $f \in H^2$, we have $\sab f=af$ and so $\|\sab\| \ge \|T_a\|=\|a\|_\infty$,
where $T_a$ is the usual Toeplitz operator defined by  $T_a f=P^+ af$.
Analogously, $\|\sab\| \ge \|b\|_\infty$ and so $\|\sab\| \ge M$. 

Next, $\|\sab\| \le   \|a\|_\infty + \|b\|_\infty$ is clear. Finally,
from
\[
\|\phi\|_2^2 = \|P^+ \phi\|^2_2 + \|P^- \phi\|_2^2, \qquad (\phi \in L^2),
\]
we know that for each such $\phi$ there is a $\theta \in [0,\pi/2]$ such that
$\|P^+ \phi\|_2 = \|\phi\|_2 \cos \theta$ and $\|P^- \phi\|_2 = \|\phi\|_2 \sin \theta$.

Thus 
\begin{eqnarray*}
\|\sab \phi\|_2 &=& \|aP^+\phi + bP^-\phi\|_2 \\
&\le& (\|a\|_\infty\cos\theta+ \|b\|_\infty \sin\theta) \|\phi\|_2 \\
&\le& \max\{ \|a\|_\infty,\|b\|_\infty\}(\cos\theta+\sin\theta) \|\phi\|_2 \le M \sqrt2 \|\phi\|_2.
\end{eqnarray*}
So we also have 
\beq\label{eq:2.5}
\|\sab\|\le \sqrt2 M,
\eeq
 and the second inequality in \eqref{eq:2.2} follows.
\endpf

As a consequence, we can characterise the zero paired operator:

\begin{cor}
$\sab=0$ if and only if $a=b=0$ in $L^\infty$.
\end{cor}

\begin{rem}{\rm
The estimates obtained in Proposition \ref{prop:2.2} are optimal. We have nontrivial cases where 
$\| \sab\|=M$, e.g., if $a=b=1$, and cases where $\| \sab\|=\sqrt 2 M$, e.g. if
$a=1$ and $b=z$. In the latter case, taking $\phi(z)=1+\dfrac{1}{z}$, we have
$\sab\phi = 2$ and $\|\phi\|=\sqrt 2$, so
$\|\sab\| \ge \sqrt 2$, and, from \eqref{eq:2.5}, $\|\sab\|=\sqrt 2=\sqrt 2 M$.
}
\end{rem}

\begin{prop}
 $\|S_{a,b}\|=\|a\|_\infty + \|b\|_\infty$
if and only if either $a$ or $b$ is $0$ a.e.
\end{prop}

\beginpf
First if $b=0$ then $S_{a,b} \phi= a P^+ \phi$ and since $\|S_{a,0}\|\ge \|T_a\|=\|a\|_\infty$
(where $T_a$ is the Toeplitz operator) it is clear that $\|S_{a,0}\|=\|a\|_\infty$ too. Similarly if $a=0$.\\

For the converse, suppose that there is a $\phi \in L^2$ with $\|\phi\|_2=1$ and
for some $\epsilon$ with $0<\epsilon<1$ we have
\[
\|aP^+ \phi + b P^- \phi\| \ge (1-\epsilon)(\|a\|_\infty + \|b\|_\infty)  .
\]
For what we are trying to prove we may assume that $a$ and $b$ are both nonzero.
There is a $\theta \in [0,\pi/2]$ such that
$\|P^+ \phi\|_2=\cos \theta$ and $\|P^-\phi\|_2=\sin \theta$.
Thus, with $A=\|a\|_\infty$ and $B=\|b\|_\infty$ we have
\[
A \cos \theta + B \sin \theta \ge (1-\epsilon) (A+B),
\]
and hence by Cauchy--Schwarz
\[
A^2+  B^2 \ge (1-\epsilon)^2  (A^2+2AB+B^2)
\]
or $(2\epsilon-\epsilon^2)(A^2+B^2) \ge 2(1-\epsilon)^2AB$.
Now suppose that $A \ge B$ (the other way round is similar). Then we have
\[
2(2\epsilon-\epsilon^2)A^2 \ge 2(1-\epsilon)^2 AB 
\]
and so
$B \le  (2\epsilon-\epsilon^2)A/(1-\epsilon)^2$, and
since we can take $\epsilon>0$ arbitrary, this means that $B=0$.
 \endpf
 
\section{Brown--Halmos type theorems}

Paired operators are closely related to Toeplitz operators in
various ways. For instance, Toeplitz operators $T_a$,
with $a \in L^\infty$, can be seen as compressions to $H^2$ of
paired operators of the form $S_{a,b}=aP^++bP^-$
or their duals $\Sigma_{a,b}=P^+a+P^-b$, with $b \in L^\infty$.
Given that, by the Brown--Halmos theorems \cite{brown-halmos} we have 
for $a, \tilde a \in L^\infty$ that
\[
T_a T_{\tilde a}=T_{a\tilde a} \iff a \in \overline{H^\infty} \hbox{  or } \tilde a \in H^\infty,
\]
one may ask, similarly, when the composition of two paired operators is a paired operator
or, more precisely, when the equations
\begin{eqnarray}
(aP^++bP^-)(\tilde aP^++\tilde bP^-) &=& a\tilde a P^+ + b \tilde b P^-
\label{eq:n3.1}
\\
 (P^+ a + P^- b)(P^+ \tilde a + P^- \tilde b) &=& P^+(a\tilde a)+P^-(b\tilde b)
 \label{eq:n3.2}
 \end{eqnarray}
 hold.
 It is clear that \eqref{eq:n3.1} and \eqref{eq:n3.2} hold if $a=b$ and $\tilde a=\tilde b$
 respectively. Otherwise we have the following (recall that nondegeneracy was defined in
 Definition \ref{def:ndg}).
 
 \begin{prop}\label{prop:n4.1}
 Let $a,\tilde a,b, \tilde b \in L^\infty$ and suppose that
 $\{a,b\}$ and $\{\tilde a,\tilde b\}$ are nondegenerate. Then
 \begin{eqnarray}
 S_{a,b}S_{\tilde a,\tilde b}=S_{a\tilde a,b\tilde b} &\iff& \tilde a\in  H^\infty \hbox{ and } \tilde b \in \overline{H^\infty}\\
 \Sigma_{a,b}\Sigma_{\tilde a,\tilde b}=\Sigma_{a\tilde a,b\tilde b}
 & \iff & a \in \overline{H^\infty}\hbox{ and } b \in H^\infty.
  \label{eq:n3.4}  
 \end{eqnarray}
 \end{prop}
 
 \beginpf
 We have that
 \[
 (aP^++bP^-)(\tilde aP^++\tilde b P^-)=a\tilde a P^+ + b\tilde b P^- + T ,
 \]
 where
 \begin{eqnarray*}
  T  &=& -aP^-\tilde a P^+ + aP^+\tilde b P^- + bP^-\tilde a P^+ - b P^+ \tilde b P^- \\
 &=&
 (a-b)(P^+\tilde b P^- - P^- \tilde a P^+).
 \end{eqnarray*}
 If $\ta \in H^\infty$ and $\tb\in \overline{H^\infty}$ then $T=0$. Conversely,
 suppose that $T=0$. Then, for all $f \in L^2$ we have $P^+\tb P^- f = P^- \ta P^+ f =0$ and
 therefore, taking $f=1$, we have $P^- \tilde a=0$, i.e., $\ta\in H^2 \cap L^\infty = H^\infty$.
 Similarly, taking $f=\bar z$, we have $P^+\tb \bar z=0$, i.e., $\tb \bar z \in \overline{H^2_0} \cap \bar z L^\infty = \bar z \overline{H^\infty}$.
 
 The equivalence in \eqref{eq:n3.4} is proved similarly, or by taking adjoints.
  \endpf
  
  It follows from Proposition \ref{prop:n4.1} that $S_{a,b}S_{\ta,\tb}=S_{\ta,\tb}S_{a,b}$
  when $a,\ta \in H^\infty$ and $b,\tb \in \overline{H^\infty}$, and
  that $\Sigma_{a,b}\Sigma_{\ta,\tb}=\Sigma_{\ta,\tb}\Sigma_{a,b}$ whenever
  $a,\ta \in \overline{H^\infty}$ and $b,\tb \in H^\infty$. More generally, we have
  the following.
  
  \begin{prop}\label{prop:A2}
  Let $a,\ta,b,\tb \in L^\infty$ with $\{a,b\},\{\ta,\tb\}$ nondegenerate.
  Then
  \[
  S_{a,b}S_{\ta,\tb}=S_{\ta,\tb}S_{a,b} \iff
  (a-b)(P^-\ta P^+ - P^+\tb P^-)
  = (\ta-\tb) ( P^- a P^+ - P^+bP^- ).
  \]
  \end{prop}
  
  \beginpf 
\[ S_{a,b}S_{\ta,\tb}-S_{\ta,\tb}S_{a,b} \]
 \begin{eqnarray*}= -aP^-\ta P^+ + aP^+\tb P^- &+& bP^-\ta P^+ - b P^+ \tb P^-\\
& +& \ta P^- a P^+ - \ta P^+ b P^- \quad -\quad  \tb P^- a P^+ + \tb P^+ b P^-
 \end{eqnarray*}
 \[ = -(a-b)(P^- \ta P^+ - P^+\tb P^-) - (\ta-\tb)(P^+bP^--P^- a P^+). \]
 \endpf
 
 \begin{cor}\label{cor:A3}
 For $\eta\in L^\infty$ and $\{a,b\}$ nondegenerate, 
 \[
 \eta S_{a,b}=S_{a,b}\eta \iff \eta \in \CC.
 \]
 \end{cor}
  \beginpf
  $\eta I= S_{\eta,\eta}$, and applying the calculation in the
  proof of Proposition \ref{prop:A2} we have
  \[
  \eta S_{a,b}=S_{a,b}\eta \iff P^- \eta P^+ = P^+ \eta P^-,
  \]
  which is true if and only if $\eta \in \CC$.
  \endpf
  
  One may then ask for which functions $f \in L^2$ do multiplication by $\eta$ and $S_{a,b}$ commute.
  
  \begin{prop}\label{prop:A4}
  For $\eta\in L^\infty$, $\{a,b\}$ nondegenerate, and $f \in L^2$,
  \begin{eqnarray*}
  \eta S_{a,b}f=S_{a,b}(\eta f) & \iff& f \in \ker  H_\eta \oplus \ker \tilde H_\eta \\
  & \iff & \eta f_+ \in H^2, \eta f_- \in  \overline{H^2_0} \\
  & \iff & \eta f_+ = P^+(\eta f) \\
  & \iff & \eta f_- = P^- (\eta f),
  \end{eqnarray*}
  where the Hankel operators $H_\eta$ and $\tilde H_\eta$ are defined by
$H_\eta=P^-\eta P^+$ and $\tilde H_\eta = P^+ \eta P^-$,
  and we write $f_\pm = P^\pm f$.
  \end{prop}
  
  \beginpf
  From Proposition \ref{prop:A2} we have
  \[
  \eta S_{a,b}f=S_{a,b}(\eta f) \iff P^-\eta f_+ = P^+ \eta f_-,
  \]
  which happens if and only if $P^-\eta f_+ = P^+ \eta f_- = 0$,
  or equivalently $\eta f_+ \in H^2$ and $\eta f_- \in \H2b$.
  
  Now $P^+(\eta f)=P^+(\eta f_+) + P^+(\eta f_-)$;
  if $\eta f_+ \in H^2$ and $\eta f_- \in \H2b$, then we have $P^+(\eta f)=\eta f_+$.
  
  Conversely, if $\eta f_+ = P^+ (\eta f)$ then $\eta f_+ \in H^2$; on the other hand, $\eta f_+ = P^+(\eta f)=P^+(\eta f_+)+P^+(\eta f_-)=\eta f_+ + P^+(\eta f_-)$,
  so $P^+(\eta f_-)=0$ and we conclude
 that $\eta f_- \in \H2b$.
  
  The last equivalence is proved analogously.
  \endpf
  
  At this point we require the theory of inner and outer functions, and
  the inner--outer factorization in $H^2$. A suitable reference is \cite[Chap. 1]{Nik1}.
  
  As an example of a class of functions $\eta \in L^\infty$ and $f \in L^2$ such that
  $P^\pm(\eta f)=\eta P^\pm f$ we have the following, where $\theta$ denotes an
  inner function and $K_\theta$ is its associated model space $K_\theta=H^2 \ominus \theta H^2$.
  
  \begin{prop}\label{prop:13apr3.5}
  If $\eta \in \overline{H^\infty}$ is such that $\eta=\alpha \bar\theta h$ with $h \in  H^\infty$
  and $\alpha$ an inner function, then for all $f \in K^\perp_\theta$ we have
  $P^\pm(\eta f)=\eta P^\pm f$.
  \end{prop}
  \beginpf
  Since $f \in K^\perp_\theta$ can be written in the form
  $f=f_- + \theta f_+$ with $f_- \in \H2b$ and $f+ \in H^2$, we have that
  \[
  \eta f=\eta(f_- + \theta f_+)=\eta f_- + \alpha \bar\theta h \theta f_+ = \eta f_- + \alpha h f_+,
  \]
  where $\eta f_- \in \H2b$ and $\alpha h f_+ \in H^2$. Therefore
$P^+ (\eta f)=\alpha h f_+ = \eta\theta f_+=\eta P^+ f$ and $P^-(\eta f)=\eta f_-= \eta P^- f$.
  \endpf
  
  \begin{rem}
  We have that $w \in L^\infty$ is a multiplier from a model space $K_{\theta_1}$ onto
  another model space $K_{\theta_2}$, i.e., $w K_{\theta_1} = K_{\theta_2}$,
  if and only if  $w \in \Hol(\DD)$ with $w=\gamma \theta_2 \overline{\theta_1}\bar w$
  where $|\gamma|=1$
 (see \cite{crofoot} and \cite{CP18a}). Thus   conjugates of the surjective multipliers
 between two model spaces provide examples of the functions $\eta$ described in Proposition
 \ref{prop:13apr3.5}.
  \end{rem}

 \section{Paired kernels}
 
 We denote the kernel of a paired operator, which we call a {\em paired kernel}, by
 \beq
 \kerab = \ker \sab.
 \eeq
Paired kernels are closely connected with Toeplitz kernels. To obtain this relation, we start
 by recalling that two operators $T$ and $S$, on the Banach spaces $X$ and $Y$ respectively,
 are said to be {\em equivalent\/} if and only if there exist invertible operators
 $E$ and $F$ such that $T=ESF$. Clearly, in this case we have the isomorphism $\ker T \simeq \ker S$ since we have
 $\ker T= F^{-1} \ker S$. The operators $T$ and $S$ are said to be 
 {\em equivalent after extension\/} if and only if there are Banach spaces $\tilde X$ and $\tilde Y$
 such that the operators
 \[
 \tilde T= \left[ \begin{matrix} T & 0 \\ 0 & I_{\tilde X} \end{matrix}\right]
 \qquad \hbox{and} \qquad
  \tilde S= \left[ \begin{matrix} S & 0 \\ 0 & I_{\tilde Y} \end{matrix}\right]
 \]
 are equivalent (see \cite{BTsk}).
 
 If $b \in \G L^\infty$ (that is, if it is invertible in $L^\infty$), then it is clear
 that $\sab$ is equivalent to $S_{G,1}$, where $G=a/b$. On the other hand,
 letting again $T_G$ denote the Toeplitz operator with symbol $G \in L^\infty$,
 \beq
 T_G: H^2 \to H^2, \qquad T_G={P^+GP^+}_{| H^2},
 \eeq
 we have that $T_G$ is equivalent after extension to the paired operator
 $S_{G,1}$ (see \cite{C}). This means in particular that the kernels of the two
 operators are isomorphic,
 \beq\label{eq:3.3}
 \ker T_G \simeq \ker S_{G,1} = \k_{G,1},
 \eeq
 the isomorphism in this case being induced by the mapping $P^+$, and we have
 \beq
 \ker T_G = P^+ \ker_{G,1} =: \k_{G,1}^+.
 \eeq
 More generally, we denote
 \beq\label{eq:3.5}
 \k^\pm_{a,b}= P^\pm \kerab = P^\pm \ker \sab.
 \eeq
 
 As happens with Toeplitz kernels, paired kernels and their projections \eqref{eq:3.5}
 possess many interesting properties which, in the case of $\k^\pm_{a,b}$, generalize those of Toeplitz kernels.\\
 
It is well known that, for $a \in L^\infty \setminus \{0\}$ we have
$\ker T_a= \{0\}$ if $a \in H^\infty$ and $\ker T_a \not= \{0\}$ if $a \in \H2b$.
 A corresponding result for paired kernels is the following.
 
 \begin{prop}
 Let $a,b \in L^\infty \setminus \{0\}$. Then:
 \begin{enumerate}
 \item[(i)] $\kerab = \{0\}$ if $a \in H^\infty$ and $b \in \overline{H^\infty}$;
 \item[(ii)] $\kerab \ne \{0\}$ if $a \in \overline{H^\infty}$ and $b \in H^\infty$ has a nontrivial inner factor;
 \item[(iii)] $\kerab \ne \{0\}$ if $a \in \overline{H^\infty_0}$ and $b \in H^\infty$.
 \end{enumerate}
 \end{prop}
 
 \beginpf
(i) $aP^+f + b P^- f = 0 \iff aP^+f=-bP^- f$ and,
since the left-hand side of this equation is in $H^2$ while the right-hand side is in $\H2b$,
both are zero and it follows that $P^\pm f=0$, so $f=0$.

(ii) Let $b=b_i b_o$ be the inner-outer factorization of $b$. Then we have
$af_+ + bf_- = 0$ for $f_+= \dfrac{b_i-b_i(0)}{z} b_o$ and $f_- = -a \dfrac{1-b_i(0)\overline{b_i}}{z}$. 

(iii) Clearly, taking $f_+=b$ and $f_-=-a$ gives a function $f=f_++f_-$ in $\kerab$.
Alternatively, the result follows from part (ii), since $\kerab=\k_{za,zb}$.
 \endpf
 
   We may see  (iii) as a generalization of the property that 
  $\ker T_a \ne \{0\}$ if $a \in \overline{H^\infty_0}$, as follows:
  
  If $a \in \overline{H^\infty_0}$, then 
  saying that $\phi_+ \in \ker T_a$ is equivalent to saying that there is a $\phi_- \in \H2b$ such that
$
  a\phi_+=\phi_- $. That is, $a\phi_+ + b \phi_- =0$, where $a \in \H2b$ and $b = -1 \in H^\infty$.\\

 We now use the notation
 \beq
 \phi_\pm = P^\pm \phi, \qquad (\phi \in L^2).
 \eeq
 \begin{rem}\label{rem:3.0}
 Note that, with the assumption that $\{a,b\}$ is nondegenerate, if $\phi\in\kerab$ and either $\phi_+=0$ or 
 $\phi_-=0$ on a set of positive measure then, by the Luzin--Privalov theorem, we must have $\phi=0$.
 \end{rem}
 Clearly, for any measurable complex-valued function $\eta$ defined a.e.\ on $\TT$ such 
 that $\eta \ne 0$ a.e.\ on $\TT$, we have the implication
 \beq\label{eq:3.7}
 \tilde a=\eta a, \quad \tilde b=\eta b\quad  \implies\quad \kerab=\k_{\tilde a,\tilde b}.
 \eeq
 It is thus natural to ask when two paired kernels are equal or related by inclusion, and whether the converse of \eqref{eq:3.7} holds.
 
 \begin{prop}\label{prop:3.1} For a nondegenerate pair $\{a,b\}$
if $\kerab \ne \{0\}$, then, for any $\tilde a,\tilde b \in L^\infty$,
it holds that
 \[
 \kerab = \k_{\tilde a,\tilde b} \quad \iff \quad a\tilde b=\tilde a b.
 \]
 \end{prop}
 \beginpf
 Suppose that $\kerab \ne \{0\}$ and $\kerab = \k_{\tilde a,\tilde b}$. Then,
 for every $\phi\in\kerab$ with $\phi \ne 0$, we have
 \beq
 a\phi_++b\phi_- =\tilde a\phi_++\tilde b \phi_- = 0,
 \eeq
 which implies that $(a\tilde b-\tilde a b)\phi_+ \phi_- = 0$.
 
 If $a\tilde b-\tilde a b \ne 0$ on a set of positive measure, then
 $\phi_+$ or $\phi_-$ must vanish on a set of positive measure, and by Remark \ref{rem:3.0} we have $\phi=0$, which is a contradiction. We
 conclude then that $a\tilde b-\tilde a b=0$.\\
 
 Conversely, suppose that $a\tilde b=\tilde a b$. Taking $\eta=\tilde a/a = \tilde b/b$ we have
 $\kerab=\k_{\tilde a,\tilde b}$ by \eqref{eq:3.7}. 
\endpf

This shows that the converse of \eqref{eq:3.7} holds. We also have the following.

\begin{cor}\label{cor:3.2}
Each nonzero $\kerab$
with $\{a,b\}$ nondegenerate is uniquely determined by any of its nonzero elements and
\[
\kerab \cap \k_{\tilde a,\tilde b} \ne \{0\} \quad \iff \quad \kerab  =  \k_{\tilde a,\tilde b}.
\]
if $\tilde a$ and $\tilde b$ are also in $L^\infty$.
\end{cor}
 \begin{cor}
 With the same assumptions as in Proposition \ref{prop:3.1},
 we have 
 \[
 \kerab \subseteq \k_{\tilde a,\tilde b} \quad \iff \quad
 \left\{
 \kerab = \k_{\tilde a,\tilde b}\quad  \vee  \quad \kerab = \{0\}
 \right\}.
 \]
 \end{cor}
 
 The following result shows that, in its turn, each
 nonzero $\phi \in L^2$ uniquely determines a paired kernel.
 
 \begin{thm}
 For each $\phi \in L^2\setminus \{0\}$, there is one and only one paired kernel to which $\phi$ belongs.
 \end{thm}
 
 \beginpf
 Take   $\phi=\phi_+ + \phi_- \in L^2 \setminus \{0\}$ and let
 $\phi_+ = I_+O_+$ and $\phi_-=I_-O_-$, where
 $I_+$ and $\overline{ I_-}$ are inner functions and $O_+$ and $\bar z \overline{O_-}$
 are outer functions in $H^2$. Then we have
 \beq\label{eq:3.9}
 \frac{I_- O_-}{I_+ O_+}\phi_+ = \phi_-; \qquad \hbox{that is,}
 \qquad \overline{I_+}\frac{\phi_+}{O_+} = \overline{I_-} \frac{\phi_-}{O_-}.
 \eeq
Now $1/O_+$ belongs to the Smirnoff class ${\mathcal N}_+$, so that there are
$H_{1+}$ and $H_{2+}$ in $H^\infty$ such that
\beq\label{eq:3.10}
\frac{1}{O_+} = \frac{H_{1+}}{H_{2+}}
\eeq
(see, for example, \cite{Nik1}). Analogously,
\beq\label{eq:3.11}
\frac{1}{\bar z \overline{O_-}}= \frac{h_{1+}}{h_{2+}}
\eeq
with $h_{1+}, h_{2+}  \in H^\infty$.
Thus, from \eqref{eq:3.9}--\eqref{eq:3.11} we have
\[
(\overline{I_+} H_{1+} \overline{h_{2+}}) \phi_+ - (\overline{I_-}   z \overline{h_{1+}} H_{2+}) \phi_- = 0
\]
so $\phi \in \kerab$ with  $a=\overline{I_+} H_{1+} \overline{h_{2+}}$ and
$b=-\overline{I_-}   z \overline{h_{1+}} H_{2+}$, which are $L^\infty$ functions.
By Corollary \ref{cor:3.2}, $\kerab$ is the only paired kernel containing $\phi$.
\endpf


\section{Invariance and near invariance}

Toeplitz kernels are $S^*$-invariant, where $S^*$ is the backward shift $S^*=T_{\bar z}$, whenever the associated symbol is in $\overline{H^\infty}$, that is
\beq
S^*(\ker T_\phi) \subseteq \ker T_\phi \quad \hbox{for} \quad \phi \in \overline{H^\infty}.
\eeq
These Toeplitz kernels are the {\em model spaces}, and they are the
nontrivial subspaces of $H^2$ that are invariant for $S^*$. 

Recall that a subspace $M \subseteq H^2$  is said to be
{\em nearly $S^*$-invariant} if 
\beq\label{eq:b2}
f \in M, \quad f(0)=0 \implies S^*f \in M.
\eeq
While not all
Toeplitz kernels are $S^*$-invariant subspaces of $H^2$, they are all
{\em nearly $S^*$-invariant\/} (see \cite{CP14} for more on this subject).
Note that for $f \in H^2$ we have $f(0)=0$ if and only if $f \in \ker H_{\bar z}$, where $H_{\bar z}$ is the Hankel operator
$H_{\bar z}=P^- \bar z P^+_{| H^2}$. Thus \eqref{eq:b2} for $M=\ker T_\phi$ can equivalently be
written as
\beq
f \in \ker T_\phi \cap \ker H_{\bar z} \implies T_{\bar z}f \in \ker T_\phi.
\eeq
This is an important property of Toeplitz kernels which can be seen as describing, in the case
when $\ker T_G$ is not invariant for $S^*$, the part of $\ker T_\phi$ that is mapped onto
$\ker T_\phi$ by $T_{\bar z}=S^*$.

More generally, Toeplitz kernels are nearly $T_\eta$-invariant for every $\eta \in \overline{H^\infty}$, i.e.,
\beq\label{eq:B4}
f \in \ker T_\phi \cap \ker H_\eta \implies T_\eta f \in \ker T_\phi
\eeq
(see \cite[Sec. 3]{CP14}), and if in addition $\phi \in \overline{H^\infty}$, then $\ker T_\phi$ is
invariant for $T_\eta$. It is natural to ask whether similar invariance properties hold for paired kernels.
We start by recalling that, for two operators $A$ and $B$,
\beq
AB=BA \implies \ker A \hbox{ is invariant for } B,
\eeq
since, if $f \in \ker A$, then $A(Bf)=B(Af)=0$ and so $Bf \in \ker A$. From this and
Proposition \ref{prop:A2} we get the following.

\begin{prop}
Let $a,\ta,b,\tb \in L^\infty$ with $\{a,b\}, \{\ta,\tb\}$ nondegenerate.
Then $\ker S_{a,b}$ is invariant for $S_{\ta,\tb}$ if $a,\ta \in H^\infty$ and
$b,\tb \in \overline{H^\infty}$; also $\ker \Sigma_{a,b}$ is invariant
for $\Sigma_{\ta,\tb}$ if $a,\ta \in \overline{H^\infty}$ and $b,\tb \in H^\infty$.
\end{prop}

Regarding the invariance properties with respect to $S_{\eta,\eta}=\eta I$
we have:

\begin{prop}
For $\eta\in L^\infty$ and $\{a,b\}$ nondegenerate we have
that if $f \in \kerab$ then we have $\eta f \in \kerab$ if and only if
\beq
f \in \ker H_\eta \oplus \ker \tilde H_\eta,
\eeq
where $H_\eta=P^-\eta P^+{}_{| H^2}$ and $\tilde H_\eta = P^+ \eta P^-{}_{| \H2b}$.
\end{prop}

\beginpf
Let $f, \eta f \in \kerab$. Then
\begin{eqnarray*}
0 &=&  aP^+\eta f + bP^-\eta f \\
&=& aP^+\eta f_+ + a P^+ \eta f_- + b P^- \eta f_+ + b P^- \eta f_- \\
&=& a\eta f_+ - a P^- \eta f_+ + aP^+\eta f_- + bP^- \eta f_+ + b\eta f_- - b P^+ \eta f_- \\
&=& -(a-b) (P^- \eta f_+ - P^+ \eta f_-),
\end{eqnarray*}
so $P^-\eta f_+ = P^+\eta f_-=0$, which is equivalent to 
$f \in \ker H_\eta \oplus \ker \tilde H_\eta$.

The converse follows from the same calculation.
\endpf

We get an analogue of \eqref{eq:B4} as follows:

\begin{prop}
For $\eta \in L^\infty$ and $\{a,b\}$ nondegenerate 
and $f \in L^2$,
\[
f \in \kerab \cap (\ker H_\eta \oplus \ker \tilde H_\eta) \implies S_{\eta,\eta} f \in \kerab.
\]
\end{prop}

\section{Coburn's lemma}

Coburn's lemma for Toeplitz operators \cite{coburn} asserts that for $g \in L^\infty$ 
at least one of $\ker T_g$ and $\ker T^*_g$
is $\{0\}$.
To obtain a generalization of this lemma    we study the relations between the kernels of
$S_{a,b}$, $S_{\bar b,\bar a}$, and $\Sigma_{\bar a,\bar b}=S^*_{a,b}$.

\begin{lem}\label{lem:3.5}
$P^\pm (\bar z \bar \phi)=\bar z \overline{\phi_{\mp}}$.
\end{lem}
\beginpf
Observe that $\bar z \bar \phi =\underbrace{ \bar z \overline{\phi_+}}_{\in \overline{H^2_0}} + 
\underbrace{\bar z \overline{\phi_-}}_{\in H^2}$
and the result follows.
\endpf

\begin{lem}
$\overline{\kerab}=z \k_{\bar b,\bar a}$.
\end{lem}
\beginpf
We have
\begin{eqnarray*}
\phi\in \kerab &\iff& a\phi_+ + b\phi_- = 0 \\
&\iff& \bar a \bar z \ol{\phi_+}+\bar b \bar z \ol{\phi_-} = 0\\
&\iff& \bar b P^+(\bar z \bar \phi) + \bar a P^-(\bar z \bar \phi) = 0\\
&\iff& \bar z \bar \phi \in \k_{\bar b,\bar a},
\end{eqnarray*}
where the third equivalence follows from Lemma \ref{lem:3.5}.
\endpf

We can thus define an antilinear isomorphism between $\kerab$ and $\k_{\bar b,\bar a}$,
as follows.

\begin{prop}\label{prop:3.7}
The antilinear operator
\beq
J: \kerab \to \k_{\bar b,\bar a}, \qquad J\phi=\bar z \bar \phi
\eeq
is well defined and bijective.
\end{prop}

Next, we construct an isomorphism from
$\k^*_{a,b}:= \ker \sab^*$ onto
$\ker_{\bar a,\bar b}$.

\begin{prop}\label{prop:apr5}
For nondegenerate $\{a,b\}$
the operator $\tilde J: \k^*_{a,b} \to \ker_{\bar a,\bar b}$ defined by
\beq
\tilde J \psi = (\bar a-\bar b)\psi, \qquad (\psi \in \k^*_{a,b})
\eeq
is injective, and if 
\beq\label{eq:3.14}
a \in \G L^\infty, \hbox{ or } b \in \G L^\infty , \hbox{ or } (a-b) \in \G L^\infty,
\eeq
then $\tilde J$ is an isomorphism from $\k^*_{a,b}$ onto $\ker_{\bar a,\bar b}$.
\end{prop}

\beginpf
We begin by showing that the operator is well defined, i.e.,
\beq
(\bar a-\bar b) \ker \sab^* \subseteq \ker S_{\bar a,\bar b}.
\eeq
Let $\psi \in  \kerab^* = \ker \sab^*$, which means that
\beq\label{eq:3.16}
P^+(\bar a \psi) = -P^- (\bar b \psi)=0,
\eeq
since $H^2 \cap \overline{H^2_0} = \{0\}$. Then
\begin{eqnarray*}
S_{\bar a,\bar b}(\bar a-\bar b)\psi &=& \bar a P^+(\bar a \psi)-\bar a P^+(\bar b \psi)+  \bar b P^-(\bar a \psi)-\bar b P^-(\bar b \psi)\\
&=& -\bar a P^+(\bar b \psi)+ \bar b P^-(\bar a \psi)\\
&=& -\bar a(I-P^-)(\bar b \psi) + \bar b (I-P^+)(\bar a \psi)\\
&=& -\bar a\bar b \psi + \bar b\bar a\psi = 0,
\end{eqnarray*}
so $(\bar a-\bar b)\psi \in \k_{\bar a,\bar b}$.

Now the operator $\tilde J$ is injective because we assume that $a-b \ne 0$ a.e.\ on $\TT$, so it is left to show that
it is also surjective when \eqref{eq:3.14} holds.

This is obvious if $a-b \in \G L^\infty$; suppose now
that $a \in \G L^\infty$ or $b \in \G L^\infty$. We have, for 
$\psi \in \kerab^*$, which must satisfy \eqref{eq:3.16},
\[
(\bar a-\bar b)\psi = \phi \iff \underbrace{\bar a \psi}_{\in \overline{H^2_0}} - 
\underbrace{\bar b \psi}_{\in H^2} =  \phi \iff \bar a \psi = \phi_- \hbox{ and } \bar b \psi = - \phi_+.
\]
Thus we have, since $\bar a \phi_+ + \bar b \phi_-=0$,
\[
\psi = \frac{\phi_-}{\bar a} = \tilde J^{-1}\phi \qquad \hbox{if} \quad a \in \G L^\infty
\]
and
\[
\psi = -\frac{\phi_+}{\bar b} = \tilde J^{-1}\phi \qquad \hbox{if} \quad b \in \G L^\infty.
\]
\endpf


\begin{cor}
With the same assumptions as in Proposition \ref{prop:apr5}, if \eqref{eq:3.14} holds, then
\[
\tilde J^{-1}: \k_{\bar a,\bar b} \to \kerab^*
\]
is defined by
\[
\tilde J^{-1} \phi =
\begin{cases}\dfrac{\phi}{\bar a-\bar b} & \hbox{if} \quad a-b \in \G L^\infty,\\ \\
\dfrac{P^-\phi}{\bar a} & \hbox{if} \quad a  \in \G L^\infty,\\ \\
-\dfrac{P^+ \phi}{\bar b} & \hbox{if} \quad b  \in \G L^\infty. 
\end{cases}
\]
\end{cor}

\begin{cor}
If $\k_{\bar a,\bar b}=\{0\}$, then $\kerab^*=\{0\}$.
\end{cor}

Taking Proposition \ref{prop:3.7} into account, we also have:

\begin{cor}\label{cor:3.11}
The following are equivalent:\\
(i) $\k_{b,a}=\{0\}$;\\
(ii) $\k_{\bar a,\bar b}=\{0\}$;\\
and, if at least one of $a,b$ and $a-b$ lies in $\G L^\infty$,\\
(iii) $\kerab^*=\{0\}$.
\end{cor}

We can now formulate a generalization of Coburn's lemma for Toeplitz operators \cite{coburn}.

\begin{thm} Suppose that $a,b \ne 0$ a.e. on $\TT$. Then\\
(i) We have $\kerab=\{0\}$ or $\k_{b,a}=\{0\}$;
equivalently $\kerab=\{0\}$ or $\k_{\bar a , \bar b}=\{0\}$.\\
(ii) If at least one of $a,b$ and $a-b$ lies in $\G L^\infty$, then
$\kerab=\{0\}$ or $\kerab^*=\{0\}$.
\end{thm}

\beginpf
Suppose that there exist $\phi,\psi \ne 0$ with $\phi \in \kerab, \psi \in \kerba$. Then
clearly
\[
a\phi_+ = -b\phi_- \qquad \hbox{and} \qquad b\psi_+ = -a \psi_-.
\]
Therefore $ab \phi_+ \psi_+ = ab \phi_- \psi_-$, and, since
$a,b \ne 0$ a.e.\ on $\TT$,
we have
\[
\underbrace{\phi_+ \psi_+}_{\in H^1} = \underbrace{ \phi_- \psi_-}_{\in \overline{H^1_0}} ,
\]
so both the terms above are $0$. This implies that either $\phi_+=0$ or $\psi_+=0$
and therefore, by Remark \ref{rem:3.0}, $\phi=0$ or $\psi=0$, which
is a contradiction.
Therefore we cannot have  both $\kerab \ne \{0\}$ and $\kerba \ne \{0\}$;
by Corollary \ref{cor:3.11}, $\kerba \ne \{0\}$ is equivalent to $\k_{\bar a,\bar b} \ne \{0\}$.

Finally, (ii) follows from Corollary \ref{cor:3.11} (iii).
\endpf

Clearly the theorem does not apply to the cases $a=0$ a.e.\ or $b=0$ a.e., as is
easily verified.

We recover Coburn's lemma when $b=1$. In that case we have
$\ker (aP^+ + P^-) = \{0\}$
or $\ker (\bar a P^+ + P^-)=\{0\}$,
where $\ker (aP^+ + P_-) \simeq \ker T_a$ and 
$\ker(\bar a P^++P^-) \simeq \ker T_a^*$ (cf. \eqref{eq:3.3}).

%

\section{Final comments}

In \cite{GHJR1,GHJR2,GHJR3}, a new class of Toeplitz-like operators $T_\omega$  is presented,
whose symbols
are complex rational functions $\omega$. The domain of $T_\omega$ is
\[
\{g_+ \in H^p: \omega g_+ = f + \rho\},
\]
where $f \in L^p$ and $\rho$ is a strictly proper (i.e., zero at $\infty$) rational function with all poles on $\TT$.
Then the authors write $T_\omega g_+ = P^+ f$. The operator is densely-defined as its domain
contains all polynomials.
Its kernel is then defined to be those $g_+$ with $\omega g_+ = f_- + \rho$.

Now, in the context of paired operators, we may take a more general $\omega =\phi/p$ where $\phi \in L^\infty$
and $p$ is a polynomial with all zeros on $\TT$.
We would naturally consider the kernel associated with the  {\em bounded\/} paired
operator $S_{\phi,-p}$, since we then have, for $S_{\phi,-p}\,g=0$,
\[
\phi g_+ = p g_-,  \qquad \hbox{that is,} \qquad \omega g_+ = g_-.
\]
In general this gives a different outcome, as can be seen with the example
$\omega(z)=1/(z-1)$. Its  $T_\omega$ kernel is the space of all $g_+$ with
$g_+/(z-1)=f_-+ \rho$.

However, the kernel of $S_{1,1-z}$ is the space of all $g$ with
$g_+ + (1-z) g_- = 0$.
Thus the kernels are different, namely $\CC$ in the first case and 
$\{0\}$ in the second.

Finally, another recent paper that studies pairs of projections 
(but in the context of model spaces, giving insights into
truncated Toeplitz operators) is \cite{DPS24}.

\section*{Acknowledgements}
This work was partially
supported by FCT/Portugal through CAMGSD, IST-ID, projects UIDB/04459/2020 and UIDP/04459/2020.
The second author would like to thank the Calouste Gulbenkian Foundation for partially funding this work under the Talentos Científicos Matem\'atica 2021 scholarship.


\begin{thebibliography}{99}



\bibitem{BTsk}
H. Bart and V.\`E. Tsekanovski\u\i,   Matricial coupling and equivalence after extension. {\em Operator theory and complex analysis (Sapporo, 1991)}, 143--160, Oper. Theory Adv. Appl., 59, Birkh\"auser, Basel, 1992. 



\bibitem{brown-halmos}
A. Brown and P.R. Halmos,  
Algebraic properties of Toeplitz operators.
{\em J. Reine Angew. Math.} 213 (1963/64), 89--102.



\bibitem{C}
M.C. C\^amara, Toeplitz operators and Wiener--Hopf factorisation: an introduction. 
{\em Concr. Oper.} 4 (2017), no. 1, 130--145.


\bibitem{CKP}
M.C. C\^amara, K. Kli\'s-Garlicka, and M. Ptak,  
(Asymmetric) dual truncated Toeplitz operators. 
{\em Operator and norm inequalities and related topics}, 429--460,
Trends Math., Birkh\dq auser/Springer, Cham, 2022.




\bibitem{CP14}
M.C. C\^amara and J.R. Partington,  Near invariance and kernels of Toeplitz operators. 
{\em J. Anal. Math.} 124 (2014), 235--260.


\bibitem{CP18a}
M.C. C\^amara and J.R. Partington,   Multipliers and equivalences between Toeplitz kernels. 
{\em J. Math. Anal. Appl.} 465 (2018), no. 1, 557--570.
 
%

\bibitem{CP20}
M.C. C\^amara and J.R. Partington,  
 Scalar-type kernels for block Toeplitz operators. {\em J. Math. Anal. Appl.} 489 (2020), no. 1, 124111, 25 pp.




\bibitem{coburn}
L. Coburn, Weyl's theorem for nonnormal operators. {\em Michigan Math. J.} 13 (1966),
285--288.

\bibitem{crofoot}
R.B. Crofoot,   Multipliers between invariant subspaces of the backward shift. 
{\em Pacific J. Math.} 166 (1994), no. 2, 225--246. 

\bibitem{DPS24}
R. Debnath, D.K. Pradhan and J. Sarkar, 
Pairs of inner projections and two applications.
{\em J. Func. Anal.} 286 (2024), 110216.


\bibitem{GHJR1}
G.J. Groenewald, S. ter Horst, J. Jaftha and A.C.M. Ran,   A Toeplitz-like operator with rational symbol having poles on the unit circle I: Fredholm properties.{\em  Operator theory, analysis and the state space approach}, 239--268, Oper. Theory Adv. Appl., 271, Birkh\"auser/Springer, Cham, 2018. 

\bibitem{GHJR2}
G.J. Groenewald, S. ter Horst, J. Jaftha and A.C.M. Ran,   
A Toeplitz-like operator with rational symbol having poles on the unit circle II: the spectrum. {\em Interpolation and realization theory with applications to control theory}, 133--154, Oper. Theory Adv. Appl., 272, Birkh\"auser/Springer, Cham, 2019.

\bibitem{GHJR3}
G.J. Groenewald, S. ter Horst, J. Jaftha and A.C.M. Ran,   
A Toeplitz-like operator with rational symbol having poles on the unit circle III: The adjoint. {\em Integral Equations Operator Theory\/} 91 (2019), no. 5, Paper No. 43, 23 pp.

\bibitem{KS23}
O. Karlovych and E. Shargorodsky,
The Coburn lemma and the
Hartman--Wintner--Simonenko theorem
for Toeplitz operators on abstract Hardy
spaces, {\em   Integral Equations Operator Theory\/} 95 (2023), no. 1, Paper No. 6, 17 pp.


\bibitem{MP}
S. Mikhlin and S. Prossdorf, {\em Singular integral operators}. Springer-Verlag, Berlin, 1986.
Translated from German by Albrecht B\"ottcher and Reinhard Lehmann, 1986.




\bibitem{Nik1}
N.K. Nikolski,  {\em Operators, functions, and systems: an easy reading. Vol. 1. Hardy, Hankel, and Toeplitz.}  Mathematical Surveys and Monographs, 92. American Mathematical Society, Providence, RI, 2002.

\bibitem{speck}
F.-O. Speck,  Paired operators in asymmetric space setting. {\em Large truncated Toeplitz matrices, Toeplitz operators, and related topics}, 681--702, Oper. Theory Adv. Appl., 259, Birkhäuser/Springer, Cham, 2017. 



\end{thebibliography}
\end{document}